\documentclass[reqno, 12pt]{amsart}
\usepackage{amsfonts}
\usepackage{amssymb,amsmath}
\usepackage{amsthm}
\usepackage{upref}

\makeatletter
\def\LaTeX{\leavevmode L\raise.42ex
    \hbox{\kern-.3em\size{\sf@size}{0pt}\selectfont A}\kern-.15em\TeX}
 \makeatother

\numberwithin{equation}{section}

\newtheorem{lemma}{Lemma}[section]
\newtheorem{theorem}[lemma]{Theorem} 
\newtheorem{corollary}[lemma]{Corollary}
\newtheorem{proposition}[lemma]{Proposition}

\theoremstyle{definition}

\newtheorem{definition}[lemma]{Definition}

\newtheorem{assumption}[lemma]{Assumption}

\newtheorem{remark}[lemma]{Remark}

 \newcommand{\supp}{\operatorname{supp}}
  \newcommand{\e}{\eqref}

\newcommand{\q}{\quad}

\newcommand{\ti}{\tilde}
\newcommand{\wt}{\widetilde}

\newcommand{\la}{\langle}
\newcommand{\ra}{\rangle}

\newcommand{\ov}{\overline}
 \renewcommand{\d}{\delta}

 \newcommand{\Ker}{\operatorname{Ker}}

\renewcommand\Im{\operatorname{Im}}
\renewcommand\Re{\operatorname{Re}}

\newenvironment{pf}{\begin{proof}}{\end{proof}}

\def\qqq{\mathrel{\subset\mkern-15mu\lower.38ex\hbox{${\scriptscriptstyle\rightarrow}$}}}

\let\cal\mathcal

\let\Bbb\mathbb

  \DeclareMathOperator{\spec}{spec}
  
\begin{document}

\title {Diagonalizations of two classes of unbounded Hankel operators}
\author{ D. R. Yafaev}
\address{ IRMAR, Universit\'{e} de Rennes I\\ Campus de
  Beaulieu, 35042 Rennes Cedex, FRANCE}
\email{yafaev@univ-rennes1.fr}


\begin{abstract}
We show that every Hankel operator $H$ is unitarily equivalent to a pseudo-differential operator $A$ of a special structure acting in   the space $L^2 ({\Bbb R}) $. As an example, 
we   consider integral operators $H$ in the space $L^2 ({\Bbb R}_{+}) $  with  kernels $P (\ln ( t+s)) ( t+s)^{-1}$ where $P(x)$ is an arbitrary real polynomial   of degree $K$. In this case,  $A$ is  a differential operator of the same order $K$. This allows us to study spectral properties of   Hankel operators $H$ with such kernels. In particular, we show that the essential
  spectrum of $H$ coincides with the whole axis for $K$ odd, and it coincides with the positive half-axis for $K$ even. In the latter case we additionally find
    necessary and sufficient conditions for the positivity of $H$. We also consider Hankel operators whose kernels have  a strong
singularity at some positive point. We show that spectra of such operators consist of the zero eigenvalue of infinite multiplicity and eigenvalues accumulating to $+\infty$ and $-\infty$. We find the asymptotics of these eigenvalues.
   \end{abstract}

\maketitle


\section{Introduction}  

{\bf 1.1.}
Hankel operators can be defined as integral operators
\begin{equation}
(H f)(t) = \int_{0}^\infty h(t+s) f(s)ds 
\label{eq:H1}\end{equation}
in the space $L^2 ({\Bbb R}_{+}) $ with kernels $h$ that depend  on the sum of variables only. We refer to the books \cite{Pe,Po} for basic information on Hankel operators. Of course $H$ is symmetric if $  h(t)=\ov{h(t)}$. There are very few cases when Hankel operators can be explicitly diagonalized. The most simple and important case $h(t)=t^{-1}$ was considered by T.~Carleman in \cite{Ca}.


Here we study a class of Hankel operators generalizing the Carleman operator. The corresponding kernels are given by the  formula
 \begin{equation}
h(t)= P (\ln t)t^{-1}
\label{eq:LOG}\end{equation} 
where  $P (x)$ is an arbitrary polynomial.
 Hankel operators $H$ with such kernels are not   bounded unless $P(x)={\rm const}$, but, for real $P(x)$,   they can be uniquely defined as self-adjoint operators. We show that the Hankel operator with kernel \e{eq:LOG} is unitarily equivalent to the differential operator
  \begin{equation}
A= v Q (D) v, \q D=i d /d \xi,
\label{eq:LOGz}\end{equation}
 in the space $L^2 ({\Bbb R}) $.
Here $v$ is the operator of multiplication by the universal function 
  \begin{equation}
v (\xi)=\frac{\sqrt{\pi}} {\sqrt{\cosh (\pi\xi)}}
\label{eq:LOGz1}\end{equation}
and the polynomial $Q(x)$ is determined by $P (x)$. The polynomials $P (x)$ and $Q (x)$ have the same degree,   and their coefficients   are linked by an explicit formula (see subs.~3.2). In particular, $Q(x)=1$ if $P(x)=1$ which yields the familiar diagonalization of the Carleman operator.

Thus    the spectral analysis of Hankel operators with kernels 
\e{eq:LOG}  reduces to  the spectral analysis   of differential operators which in principle is very well developed. However   operators \e{eq:LOGz} are somewhat unusual because the function $ v (\xi)$ tends to zero exponentially as $|\xi|\to \infty$ so that there is a strong degeneracy  at infinity. Nevertheless we  describe completely the essential spectrum of differential operators \e{eq:LOGz}  under rather general assumptions on the function $v(\xi)$. We show that $\spec_{\rm ess} (A)={\Bbb R}$ if $K:=\deg P$ is odd, and $\spec_{\rm ess} (A)= [0, \infty)$ if $K $ is even. Moreover, it turns out that zero  is never an eigenvalue of $A$. In the case of even $K$ we also find necessary and sufficient conditions for the positivity of $A$ and for the infinitude of its negative spectrum.   Of course the same spectral results are true for   Hankel  operators $H$ with kernels  \e{eq:LOG}. For real polynomials $P(x)$ of first order, our approach yields the explicit diagonalization of Hankel operators $H$. In particular, we show that in this case the spectrum of $H$ is absolutely continuous, has multiplicity $1$ and covers the whole real line.


Actually, the unitary equivalence of the operators $H$ and $A$ is quite explicit. Let $M: L^2 ({\Bbb R}_{+}) \to L^2 ({\Bbb R} )$ be the Mellin transform; it is a unitary mapping. Set
  \begin{equation}
(F f) (\xi)= \frac{\Gamma(1/2+i\xi)}{|\Gamma(1/2+i\xi)|} (M f)(\xi)
\label{eq:MAID1}\end{equation}
where $\Gamma(\cdot)$ is the gamma function.
We show that
  \begin{equation}
H= F^\ast A F.
\label{eq:LOGz3}\end{equation}

Our proof of this identity follows the approach of 
 \cite{Y} where general Hankel operators $H$ were considered. For an arbitrary $H$, the function $Q(x)$ in formula \e{eq:LOGz} is a distribution which may be  (for example, for finite rank $H$) very singular. However it is a polynomial for Hankel operators with kernels \e{eq:LOG} so that $A$  is the explicit differential operator in this case.
 
    \medskip
    
     {\bf 1.2.}
  Kernels \e{eq:LOG}  are singular at the points $t=0$ and $t=\infty$.  We also consider another class of kernels which are singular at some point $t_{0} >0$.  
   We assume that 
   \begin{equation}
h(t)=\sum_{k=0}^K h_{k} \d^{(k)} (t-t_{0}), \q h_{k}=\bar{h}_{k}, \q h_{K}\neq 0, 
\label{eq:DD}\end{equation}
where  $\d(\cdot)$ is the Dirac delta function. It turns out that Hankel operators with such kernels  reduce to ``differential" operators with the reflection and a shift of the argument.

Spectral properties of Hankel operators with kernels \e{eq:LOG}  and \e{eq:DD} are completely different.  As discussed in \cite{Y},  Hankel operators can be sign-definite only for   $h\in C^\infty ({\Bbb R}_{+})$  which is of course not true for kernels \e{eq:DD}. If $K=0$, then, as shown in \cite{Y}, the spectrum of $H$ consists of three eigenvalues $0$, $h_{0}$ and $-h_{0}$ of infinite multiplicity each. 
We shall prove here that for $K\geq 1$   the spectrum of the operator $H$ with kernel  \e{eq:DD} consists of the zero eigenvalue of infinite multiplicity and an infinite number of  
eigenvalues of finite multiplicities accumulating both at $+\infty$ and $-\infty$. Moreover, we shall find the leading term of asymptotics of these eigenvalues. 

Recall that a symbol of a Hankel operator with kernel $h(t)$ is defined as a function $\omega(\lambda)$,  $\lambda\in {\Bbb R}$, such that $(2\pi)^{-1/2} (\Phi\omega) (t)=h(t)$ where   $\Phi$  is the Fourier transform. 
 Since by the Nehari theorem, Hankel operators with bounded symbols are bounded, the symbols of operators with kernels \e{eq:H1} and \e{eq:DD} are necessarily unbounded functions.  For kernels \e{eq:H1} symbols can be constructed by the formula
 \[
 \omega(\lambda)=2i \int_{0}^\infty \sin (\lambda t) h(t)dt= 2i \int_{0}^\infty  \sin t  \,   P(\ln t -\ln\lambda)\, t^{-1} dt
 \]
 for $\lambda>0$ and $ \omega(-\lambda)=- \omega(\lambda)$. Thus $ \omega(\lambda)$ is a $C^\infty$ function for $\lambda\neq 0$ with logarithmic singularities at $\lambda=0$ and $\lambda=\infty$. For kernels \e{eq:DD}, the symbol equals
 \[
 \omega(\lambda)=\sum_{k=0}^K h_{k} (-i\lambda)^k e^{i\lambda t_{0}},
 \]
 so that it is a  $C^\infty$ function with a power growth and an oscillation as $|\lambda|\to \infty$.

    \medskip

 {\bf 1.3.}
  Let us introduce some standard
 notation. 
 We denote by $\Phi$,
 \[
(\Phi u) (\xi)=  (2\pi)^{-1/2} \int_{-\infty}^\infty u(x) e^{ -i x \xi} dx,
\]
 the Fourier transform.
     The space $\cal Z= \cal Z ({\Bbb R})$ of test functions  is defined as the subset  
 of the Schwartz  space ${\cal S}={\cal S} ({\Bbb R}) $ which consists of functions $u $ admitting the analytic continuation to   entire functions in the   complex plane $\Bbb C$   and satisfying   bounds 
 \[
  | u (z)| \leq C_{n}  (1+| z |)^{-n} e^{r |\Im z |}, \q z\in \Bbb C,
  \]
  for some $r=r(u)>0$  and all $n$. We recall (see, e.g., the book \cite{GUECH}) that the Fourier transform 
  $\Phi : \cal Z \to C_{0}^\infty ({\Bbb R})$
  and
  $\Phi^*: C_{0}^\infty ({\Bbb R}) \to \cal Z$. 
     The dual classes of distributions (continuous antilinear functionals on ${\cal S}$, $C_{0}^\infty ({\Bbb R})$ and ${\cal Z}$) are denoted ${\cal S}'$, $C_{0}^\infty ({\Bbb R})'$ and ${\cal Z}'$, respectively. 
  We use the notation   ${\pmb\la} \cdot, \cdot {\pmb\ra}$ and $\la \cdot, \cdot\ra$ for   the  
  duality symbols in $L^2 ({\Bbb R}_{+})$ and $L^2 ({\Bbb R})$, respectively. They are   linear in the first argument and antilinear in the second argument.

     We denote by   $ {\sf H}^K({\cal J})$  the Sobolev space of functions defined on an interval ${\cal J}\subset {\Bbb R}$; 
     $ C_{0}^K({\cal J})$ is the class of $k$-times continuously differentiable functions with compact supports in ${\cal J}$.
       We often use the same notation for a function and the operator of multiplication by this function.
   The letters $c$ and $C$ (sometimes with indices) denote various positive constants whose precise values are inessential.
     


   Let us briefly describe the structure of the paper. We collect necessary results of \cite{Y}   in Section~2.   In Section~3 we establish the unitary equivalence of the Hankel operator $H$ with kernel \e{eq:LOG} and differential operator \e{eq:LOGz}. Spectral properties of the operators $A$ and hence of $H$ are studied in Section~4. We emphasize that our results on the operator $A$ do not require specific assumption \e{eq:LOGz1}. Finally,  Hankel operators $H$ with kernels \e{eq:DD} are studied in Section~5. Our presentation in this section is independent of general results of  Section~2. On the other hand, it is rather similar to that in  Section~6 of \cite{Y} where Hankel operators with discontinuous kernels  (but not as singular as  kernels  \e{eq:DD}) were considered.

\section{Hankel and pseudo-differential  operators}  

In this section  we show that an arbitrary Hankel    operator $H$ is unitarily equivalent to a pseudo-differential  operator $A$ defined by formula \e{eq:LOGz} with a distribution $Q(x)$. Our presentation   is close to  \cite{Y}, but we here  insist  upon the unitary equivalence of the operators $H$ and $A$.

\medskip

{\bf 2.1.}
Let us consider a Hankel operator $H$ defined by equality \e{eq:H1} in the space $L^2 ({\Bbb R}_{+})$.   Actually, it is more convenient to work with sesquilinear forms  instead of operators. Let us introduce the   Laplace convolution
\[
( \bar{f}_{1}\star f_{2})(t)=
\int_{0}^t    \overline{f_{1}(s)} f_{2}(t-s) ds  
\]
 of   functions $ \bar{f}_{1}$ and $ f_{2}$. Then
\begin{equation}
(Hf_{1} ,f_{2})=  {\pmb \la} h, \bar{f}_{1}\star f_{2} {\pmb \ra} 
 \label{eq:HH}\end{equation}
 where we write ${\pmb \la} \cdot, \cdot {\pmb \ra}$ instead of $( \cdot, \cdot )$ because $h$ may be a distribution.
 
We consider form \e{eq:HH} on elements ${f}_{1},  f_{2} \in \cal D$ where the set $\cal D$ is defined as follows. Put  
\[
(U f)(x) =e^{x/2} f(e^x).
\]
Then    $U: L^2 ({\Bbb R}_{+})\to L^2 ({\Bbb R} )$  is the unitary operator, and
  $f \in \cal D$ if and only if $  Uf  \in {\cal Z}$. Since
$f(t)= t^{-1/2}  (Uf) (\ln t)$
and ${\cal Z}\subset {\cal S}$, we see that   functions $f  \in \cal D$ and their derivatives satisfy the estimates
  \begin{equation}
|f^{(m)}(t)|=  C_{n,m} t^{-1/2-m} (1+|\ln t|)^{-n}
  \label{eq:assf}\end{equation}
for all $n$ and $m$. Of course,  the set $\cal D$  is dense in the space $L^2 ({\Bbb R}_{+})$.
It is shown in \cite{Y} that  if  ${f}_{1},  f_{2} \in \cal D$, then  the function
 \[
\Omega(x)=(\bar{f}_{1}\star f_{2} )(e^x) 
\]
 belongs to the space  ${\cal Z}$.

  With respect to $h$, we assume that the distribution  
 \begin{equation}
 \theta (x)= e^x h(e^x)
 \label{eq:HH3}\end{equation}
 is an element of the space $\cal Z'$. The set of all such  $h$ will be denoted   $ {\cal Z}'_{+}$,   that is,
 \[
h\in  {\cal Z}_{+}' \Longleftrightarrow \theta\in  {\cal Z}'.
\]
It is shown in \cite{Y} that {\it this condition is satisfied for all bounded Hankel operators }$H$. Since $\Omega\in {\cal Z}$, the form
 \[
  {\pmb\la} h, \bar{f}_{1}\star f_{2}  {\pmb\ra} = \int_{0}^\infty h( t ) ( {f}_{1}\star \bar{f}_{2} )(t) d t
  = \int_{-\infty}^\infty \theta (x ) \overline{\Omega(x)} d x=: { \la} \theta, \Omega { \ra}
\]
 is correctly defined for all ${f}_{1},  f_{2} \in \cal D$.

  Note that $h\in  {\cal Z}_{+}'$ if $h\in L^1_{\rm loc}({\Bbb R}_{+})$ and the integral
 \begin{equation}
\int_{0}^\infty | h(t)| (1+| \ln t |)^{-\kappa } dt< \infty
 \label{eq:ass}\end{equation}
 converges  for some $\kappa $. In this case the corresponding function \e{eq:HH3} satisfies the condition
  \[
\int_{-\infty}^\infty |\theta (x)| (1+| x |)^{-\kappa } dx< \infty,
 \]
 and hence $\theta \in{\cal S}'\subset {\cal Z}'$.

    \medskip
 
 {\bf 2.2.}
 Let us now give the definitions of the $b$- and $s$-functions   of a Hankel operator $H$. We formally define
  \begin{equation}
 b(\xi) = \frac{1}{ 2\pi }  \frac{ \int_{0}^\infty   h (t) t^{ -i \xi}   dt}{ \int_{0}^\infty   e^{-t} t^{ -i \xi}   dt} .
\label{eq:Bb}\end{equation}
Of course $  b(-\xi)=\ov{ b(\xi)}$ if $  h(t)=\ov{ h(t)}$. 
We call $b(\xi)$ the $b$-function of a Hankel operator $H$. Formula \e{eq:Bb} can be rewritten as
\begin{equation}
 b(\xi)  =  (2\pi)^{-1/2} a(\xi)  \Gamma (1-i\xi)^{-1}
\label{eq:M9}\end{equation}
where 
  \begin{equation}
 a(\xi) =  (\Phi \theta)(\xi)=  (2\pi)^{-1/2}\int_{0}^\infty   h (t) t^{ -i\xi}   dt 
\label{eq:M6}\end{equation}
is the Fourier transform of function  \e{eq:HH3}.

Recall that the gamma function $\Gamma(z)$ is a holomorphic function in the right half-plane and  $\Gamma(z)\neq 0$ for all $z\in{\Bbb C}$.    
 According to  the Stirling formula  the gamma function $\Gamma(z)$ tends to zero exponentially as $|z|\to\infty$ parallel with the imaginary axis. To be more precise,
    we have
\[
   \Gamma( \alpha+i \xi ) = e^{\pi i (2\alpha-1)/4} (2\pi / e)^{1/2}  \xi^{\alpha-1/2} e^{i\xi (\ln \xi  -1)} e^{-\pi \xi /2}\big(1+O(\xi^{-1})\big)  
\]
for a fixed $\alpha>0$ and $ \xi \to +\infty$. We also note that $ \Gamma(\alpha - i \xi)=\overline{\Gamma(\alpha +i \xi)}$ and
   \[
 |\Gamma(1/2+i \xi)|^2= \frac{\pi}{ \cosh(\pi\xi)}.
\]

Since the denominator in \e{eq:M9} tends to zero exponentially as $|\xi|\to\infty$, $b(\xi) $ is a ``nice" function only under very stringent assumptions on $a(\xi)$ and hence on $h(t)$. Therefore we are obliged to work with distributions which turn out to be very convenient. The Schwartz class is too restrictive for our purposes because of the exponential decay of $\Gamma (1-i\xi)$. Therefore we assume that $a\in C_{0}^\infty ({\Bbb R})'$; in this case $b(\xi)$ belongs to the same class. Our assumption on $a$ means that $\theta \in \cal Z'$ or equivalently $h \in \cal Z_{+}'$.

  Thus we are led to the following

\begin{definition}\label{HBS}
Let $h \in  {\cal Z}'_{+}$. The distribution $b\in C_{0}^\infty ({\Bbb R})'$   defined  by formulas \e{eq:HH3}, \e{eq:M9} and \e{eq:M6} is called the $b$-function of the Hankel operator $H$  (or   of its kernel $h(t)$). Its  Fourier transform $s= \sqrt{2\pi}    \Phi^* b\in {\cal Z} '$ is called the $s$-function or the {\it  sign-function} of $H$.
  \end{definition}

 Let the unitary mapping  $F : L^2 ({\Bbb R}_{+}) \to  L^2 ({\Bbb R})$ be defined by formula \e{eq:MAID1} where
  $ M= \Phi U $ is the Mellin transform. If $f \in  {\cal D}$, then $Uf\in{\cal Z}$ and hence the function $F f \in C_{0}^\infty ({\Bbb R})$. We recall that
   the function $v(\xi)$ was defined by formula \e{eq:LOGz1} and set $({\cal J}g)(\xi)=g(-\xi)$.

The following result was obtained in \cite{Y}. 

 \begin{theorem}\label{1}
 Suppose that  $h \in  {\cal Z}'_{+}$, and let  $b\in C_{0}^\infty ({\Bbb R})'$   be the corresponding $b$-function.   
  Let $f_{j}\in {\cal D}$, $j=1,2$. Then    $ g_{j} = F f_{j} \in C_{0}^\infty ({\Bbb R})$ and the representation 
\begin{equation}
{ \pmb\la} h ,  \bar{f}_{1}\star f_{2} {\pmb \ra} =    \la b,   (v {\cal J}\bar{g}_1) * (v g_{2}) \ra
\label{eq:M15s}\end{equation}
  holds.
 \end{theorem}

  Passing in the right-hand side of \e{eq:M15s} to the Fourier transforms and using that
 \[
 \Phi^*\big(   (v {\cal J}\bar{g}_1) * (v g_{2}) \big)=  (2\pi)^{1/2} \ov{\Phi^* ( v g_1)} \Phi^* (v g_2),    
 \]
 we obtain
 
   \begin{corollary}\label{four}
 Let $s \in {\cal Z} '$ be the sign-function of $h$, and let  $u_{j}= \Phi^* ( v F f_{j} ) \in {\cal Z} $. Then
 \[
{ \pmb\la} h ,  \bar{f}_{1}\star f_{2} {\pmb \ra} =   \la s,   \bar{u}_{1} u_2\ra   . 
\]
 \end{corollary}

 We note that, formally, the identity \e{eq:M15s} can be rewritten as relation \e{eq:LOGz3}
   where $A$ is the ``integral operator" with kernel
 $ v(\xi) b(\xi-\eta) v(\eta)$.
 To put it differently,   
 \begin{equation}
A= v \Phi s \Phi^* v,
\label{eq:MM2}\end{equation}
that is, $A$ is the  pseudo-differential  operator  defined by the amplitude
$ v(\xi) s(x) v (\eta)$.   We emphasize that in general $s(x)$ is a distribution so that formula \e{eq:MM2} has only a formal meaning.      According to relation \e{eq:LOGz3} a study of the operator $H$ reduces to  that of the operator $A$.

       \medskip
 
 {\bf 2.3.}
       For an arbitrary  distribution  $h \in  {\cal Z}_{+}'$, we have constructed in the previous subsection its sign-function $s  \in  {\cal Z}'$. 
It turns out that, conversely, the kernel  $h(t)$ can be recovered
from  its sign-function $s(x)$. It is convenient to introduce the distribution  
 \begin{equation}
 h^\natural (\lambda)= \lambda^{-1} s (-\ln \lambda).
 \label{eq:B}\end{equation} 
 Note that the inclusions  $ s\in {\cal Z}'$ and $  h^\natural\in {\cal Z}_{+}'$ are equivalent.
 The proof of the following result can be found in \cite{Y}. 

    \begin{theorem}\label{round}
  Let    $h \in  {\cal Z}_{+}'$, and let     $ s\in {\cal Z}'$  be the corresponding sign-function 
  $($see Definition~\ref{HBS}$)$.   Then  $h$ can be recovered from function \e{eq:B}  by the  formula 
   \begin{equation}
  h (t)  =    \int_{0}^\infty     e^{-t \lambda} h^\natural (\lambda) \lambda d \lambda.
   \label{eq:conv1}\end{equation}
   \end{theorem}
       
      Formula \e{eq:conv1} is understood of course in the sense of distributions.     We emphasize         the mappings $h\mapsto h^\natural $ as well as its inverse  $ h^\natural \mapsto h$ are one-to-one continuous mappings of the space $ {\cal Z}_{+}'$ onto itself.


   
\section{Quasi-Carleman and differential  operators}  
 
 {\bf 3.1.}
 Now we are in a position to consider  Hankel operators    with kernels 
  defined by formula
  \e{eq:LOG} where
 \begin{equation}
P (x)=\sum_{k=0}^K p_{k} x^k, \q p_{K}\neq 0,
\label{eq:LOG4}\end{equation}
is a polynomial with   coefficients $p_{k}$, $k=0,1,\ldots, K$. It is easy to see that such operators (they will be denoted by $H_0$) are well defined on the set  ${\cal D}_{0}$ of functions $f(t)$ satisfying estimate \e{eq:assf} for $m=0$ and some  $n> K+1$. Indeed,
by the Schwarz inequality for an arbitrary $\varepsilon>0$, we have the estimate
 \begin{multline}
\int_{0}^\infty dt \big|\int_{0}^\infty ds \frac{|\ln (t+s)|^k}{t+s}s^{-1/2} (1+|\ln s|)^{-n}\big|^2
\\
\leq C \int_{0}^\infty dt \int_{0}^\infty ds \frac{|\ln (t+s)|^{2k}}{(t+s)^2}  (1+|\ln s|)^{-2n+1+\varepsilon} .
  \label{eq:LOG4w}\end{multline} 
  Let us make the change of variables $(t,s)\mapsto (\tau,s)= (t+s,s)$ in the right-hand side   and   integrate first over $\tau\geq s$. Then using inequality 
  \[
  \int_s^\infty     |\ln \tau|^{  2k }\tau ^{-2}d\tau \leq C_{1} 
   (1+|\ln s|)^{  2k } s^{-1},
  \]
  we see that  expression \e{eq:LOG4w} is bounded by the integral
  \[
   C_{2}  \int_{0}^\infty     (1+|\ln s|)^{-2n + 2k+1+\varepsilon}s^{-1}ds
  \]
 which converges if $n> k+1+\varepsilon/2$.
It follows that $H_0 f\in L^2 ({\Bbb R})$ for $f\in {\cal D}_{0}$. Moreover, using the Fubini theorem, we see that $(H_0 f_{1}, f_{2})= (f_{1}, H_0 f_{2})$ for $f_{1}, f_{2}\in {\cal D}_{0}$
if all coefficients $p_{k}$ are real. 

Let us formulate the results obtained.

 \begin{lemma}\label{1zq}
 Let the  kernel $h(t)$ of a Hankel operator $H_0$ be given  by formulas \e{eq:LOG} and \e{eq:LOG4}. Then $H_0$ is well defined on the set  ${\cal D}_{0}$, and it is symmetric on ${\cal D}_{0}$ if 
 all coefficients $p_{k}$ are real. 
 \end{lemma}

  As we shall see below, 
the operator $H_0$ is essentially self-adjoint   (see, e.g., the book \cite{BS}, for background information on the theory of self-adjoint extensions of symmetric operators). The proof of this result as well as our study of spectral properties of the closure $\bar{H}_0=:H $ of $H_0$ rely on the identity \e{eq:LOGz3}. We emphasize however that the proof of \e{eq:LOGz3} does not require the assumption $p_{k}=\bar{p}_{k}$, $k=0,1,\ldots, K$.  The symmetricity of $H_0$ on the domain ${\cal D}$ is also a consequence of \e{eq:LOGz3} so that the direct proof of Lemma~\ref{1zq} could be avoided.

 \medskip
 
  {\bf 3.2.}
  Since kernels  \e{eq:LOG} satisfy condition \e{eq:ass}  with any $\kappa > K+1$, Theorem~\ref{1} can be directly applied   in this case. We only have to calculate the corresponding $b$- and $s$-functions. If $h_{k}(t)= t^{-1} \ln^k t $, then   the  function \e{eq:HH3}  equals $\theta_{k}(x)=x^k$ and its Fourier transform equals
  \[
a_{k}(\xi)= (\Phi \theta_{k}) (\xi)=
(2\pi)^{1/2} i^k \d^{(k)} (\xi).
\]
To simplify notation, we set
$
\omega ( z)=\Gamma (1-z)^{-1}.
$
 Then   function \e{eq:M9}  equals
 \[
b_{k}(\xi)=  i^k \omega ( i \xi) \d^{(k)} (\xi)= 
  \sum_{\ell =0}^k i^\ell  C^\ell _{k} \omega^{(k-\ell )} (0) \d^{(\ell )}  (\xi)
\]
where $C^\ell_{k}$ are the binomial coefficients. 

It follows that the $b$-function of kernel \e{eq:LOG}, \e{eq:LOG4} is given by the formula
 \[
b (\xi)=    \sum_{k=0}^K q_{k}  i^k \d^{(k)}   (\xi)
\]
where
 \begin{equation}
q_{k}=    \sum_{\ell =k}^K  C^k_{\ell} \omega^{(\ell-k)} (0) p_{\ell} ,\q k=0,\ldots, K, \q \omega ( z)=\Gamma (1-z)^{-1}.
\label{eq:LOG5a}\end{equation}  
It means that the operator $A$ acts by formula \e{eq:LOGz} where  
 \begin{equation}
 Q(x)=  \sum_{ k=0}^K q_{k} x^k .
\label{eq:LO4}\end{equation}
Thus for kernels   \e{eq:LOG} the sign-function  $s(x)=Q(x)$ is the polynomial. Note that according to general formula  \e{eq:conv1}, $P (x) $ can be recovered from $Q (x) $ by the equality
    \begin{equation}
 P (\ln t) =  t \int_{0}^\infty  Q (-\ln \lambda) e^{-t\lambda} d\lambda.
\label{eq:LOG5f}\end{equation}
 
Observe  that $q_{K}= p_{K}$ for all $K$. Recall that $\Gamma'(1)= - \gamma$ (the Euler constant) and
 $
 \Gamma''(1)=\gamma^2+\pi^2/6.
 $
 Therefore we have
  \begin{equation}
q_{0}= p_{0}- \gamma p_{1}, \q {\rm if} \q K=1,
\label{eq:LOG61}\end{equation}
and
  \begin{equation}
q_{0}= p_{0}- \gamma p_{1}+ (\gamma^2-\pi^2/6)p_{2},\q q_{1}=p_{1}-2\gamma p_{2}, \q {\rm if} \q K=2.
\label{eq:LOG6}\end{equation}

The following assertion is a  particular case of Theorem~\ref{1}.

 \begin{theorem}\label{1z}
 Let a kernel $h(t)$ be defined by formulas \e{eq:LOG} and \e{eq:LOG4}. Let $Q(x)$ be   polynomial \e{eq:LO4} with coefficients \e{eq:LOG5a}, and let $A$ be   differential operator \e{eq:LOGz}. Then for all $f_{j}\in {\cal D}$, $j=1,2$, the identity 
 \begin{equation}
(H f_{1}, f_{2})=  (A F f_{1}, F f_{2})
 \label{eq:MM1}\end{equation} 
holds.
 \end{theorem}

 Note that, for $h(t)=t^{-1}$,  the identity \e{eq:MM1} yields the familiar diagonalization of the Carleman operator. Indeed, in this case we have 
      \[
\theta (x)=1, \q  a(\xi)=    (2\pi )^{1/2}  \d (\xi), \q b(\xi)=   \d (\xi), \q s(x)=1.
\]
Therefore the identity \e{eq:MM1}  reads as
 \[
 (H f_1, f_2 )= 
 \int_{-\infty}^\infty \frac{\pi} {\cosh (\pi \xi)} \ti{f}_{1}(\xi)    \overline{ \ti{f}_{2}(\xi)} d\xi
    \]
    where $\ti{f}_j= M f_{j} = \Phi U f_j$, $j=1,2$,  is the Mellin transform of $f_j $.

  We emphasize that Theorem~\ref{1z} does not require that the coefficients of $P(x)$ be real.

 \medskip
 
  {\bf 3.3.}
   In view of Theorem~\ref{1z}  spectral properties of the Hankel operator $H$ are the same as those of the differential operator $A$. Therefore we forget for a while Hankel operators and study differential operators $A$ defined by formula \e{eq:LOGz}, but we not assume that the function $v (\xi)$ has  special form \e{eq:LOGz1}. 
 {\it    We suppose that $v =\bar{v}  \in C ^K ({\Bbb R})$ and that the  coefficients of the polynomial $Q(x)$ of degree 
   $K$ are real and} $q_{K}\neq 0$.  Then the operator $A_{0}$  defined by formula \e{eq:LOGz}  on the domain $C_{0}^K ({\Bbb R})$ is symmetric in $L^2({\Bbb R})$. We emphasize that operators \e{eq:LOGz} require a special study because the function $v (\xi) $ may tend to zero as $|\xi|\to \infty$. 
   
   Let us start with the case $K=1$ when $Q(x)= q_{0}+q_{1}x$ and $A_{0}$ can be standardly reduced by a change of variables and a gauge transformation
   to the differential operator $ q_{1} D $. We suppose that  $v(\xi) > 0$ and
    \begin{equation}
\int_{{\Bbb R}_{+}}v(\xi)^{-2} d\xi= \int_{{\Bbb R}_{-}}v(\xi)^{-2} d\xi= \infty.
\label{eq:int}\end{equation} 
Under this assumption the operator $T  :L^2({\Bbb R}) \to L^2({\Bbb R})$   defined by the relation
     \begin{equation}
  ( T  g )(\xi)= v(\xi)^{-1} e^{iq_{0}q_{1}^{-1} \xi } g(\int_{0}^\xi v(\eta)^{-2} d\eta)
\label{eq:int1}\end{equation} 
is unitary and
 \[
  A_{0}T g = i q_{1}  T g' ,\q g \in C_{0}^1 ({\Bbb R}).
  \]
  Recall that $D=id/d\xi$.
   Let the set  ${\cal D}_{*}\subset L^2 ({\Bbb R})$ consist of functions $g\in {\sf H}_{\rm loc}^1({\Bbb R}) $ such that 
  $ vQ( D) (v g) \in L^2 ({\Bbb R})$. It  is easy to see that ${\cal D}_*= T {\sf H}^1 ({\Bbb R})$.
  Thus we are led to the following assertion.

   \begin{lemma}\label{ord1}
    Suppose that $v\in C^1 ({\Bbb R})$, $v(\xi)> 0$ and condition \e{eq:int} is satisfied.
 Then  the operator $A_{0}=v( q_{0}+q_{1} D) v $   is essentially self-adjoint on $C_{0}^1 ({\Bbb R})$, and its closure $\bar{A}_{0}=:A$ is self-adjoint on the domain ${\cal D}_* =: {\cal D}(A) $. 
 The spectrum of the operator $A$ is simple, absolutely continuous, and it coincides with $\Bbb R$.
 \end{lemma}

     \begin{remark}\label{ord2}
   If both integrals \e{eq:int}  are finite, then $A_{0}$ reduces to the operator $q_{1}D$ on a finite interval. Its deficiency indices equal $(1,1)$.    If only one of integrals \e{eq:int}  is finite, then $A_{0}$ reduces to the operator $q_{1}D$ on a half-axis. Its deficiency indices equal $(0,1)$ or $(1,0)$. 
 \end{remark} 
  
 As a by-product of our considerations, we obtain the following result. It is simple but perhaps was never explicitly mentioned.
 
    \begin{proposition}\label{em}
 Suppose that $v\in C^1 ({\Bbb R})$ and $v(\xi)>0$. Let the space ${\cal K}={\cal K}_{v}$ consist of functions $ g \in {\sf H}_{\rm loc}^1({\Bbb R}) $  with the norm
   \[
 \| g \|^2_{v} = \int_{- \infty}^\infty \big( v^2(\xi) |g' (\xi)|^2+ v^{-2}(\xi) |g (\xi)|^2\big) d\xi.
   \]
   Then the set $C_{0}^1 ({\Bbb R})$ is dense in ${\cal K}$ if and only if condition \e{eq:int} is satisfied.
 \end{proposition} 

 \begin{pf}
 Let us make the change of variables
 \[
 \ti{g}(\ti{\xi})= g (\xi)\q {\rm where} \q \ti{\xi}= \int_{0}^\xi v(\eta)^{-2} d\eta
 \]
 and put
 \[
 v_{\pm}= \pm \int_{{\Bbb R}_\pm}v(\xi)^{-2} d\xi.
 \]
 Then
  \[
 \| g \|^2_{v} = \int_{v_{-}}^{v_{+}} \big(   | \ti{g}'(\ti{\xi}) |^2+   | \ti{g }(\ti{\xi})|^2\big) d\ti{\xi}.
  \]
   Since $g \in C_{0}^1 ({\Bbb R})$ if and only if $ \ti{g } \in C_{0}^1 (v_{-}, v_{+})$, it remains to use that  the set $C_{0}^1 (v_{-}, v_{+})$ is dense in the Sobole space ${\sf H}^1 (v_{-}, v_{+})$ if and only if $v_{-}=-\infty$ and $ v_{+}=+\infty$.
  \end{pf}

  Returning to the Hankel operator $H$ with kernel \e{eq:LOG} and using Theorem~\ref{1z} and equality \e{eq:LOG61}, we obtain the following result.
  
   \begin{theorem}\label{SpB11}
 Suppose that $P(x)=p_{0}+ p_{1} x$ where $p_{0}= \bar{p}_{0}$ and $p_{1} = \bar{p}_1\neq 0$. Then
 \[
 H=q_{1} F^* T   D T^* F
 \]
 where $T $ is  defined by formula \e{eq:int1}  with $q_{0}=  p_{0}-\gamma p_{1}$, $q_{1}=p_{1}$ and $v(\xi)$ is function  \e{eq:LOGz1}. The operator $H$ is essentially self-adjoint on the set ${\cal D}_{0}$,  and it is self-adjoint on the set $F^* {\cal D}_*$.  The spectrum of the operator $H$ is simple, absolutely continuous, and it coincides with $\Bbb R$.
 \end{theorem}

 \medskip
 
  {\bf 3.4.}
  Let us   pass to the case $K\geq 2$.  We recall that the operator $A_{0}= v Q(D) v$ is symmetric in $L^2 ({\Bbb R})$ on $C_{0}^K ({\Bbb R})$. Let us use the notation ${\cal A}_{0}$ for the same operator considered as a mapping ${\cal A}_{0}: C_{0}^K ({\Bbb R}) \to L^2 ({\Bbb R})$. The operator ${\cal A}_{0}^*:     L^2 ({\Bbb R}) \to C_{0}^K ({\Bbb R})'$ is defined by the relation
  \begin{equation}
 (  {\cal A}_{0} g,y )  =  \la   g , {\cal A}_{0}^* y \ra, \q g \in C_{0}^K ({\Bbb R}),\;  y \in L^2 ({\Bbb R}),
 \label{eq:ad}\end{equation}
  and is given by the same differential expression \e{eq:LOGz} where derivatives are   understood in the sense of distributions.
  
  It is also quite easy to construct the   operator $A_{0}^*$ adjoint to $A_{0}$ in the space  $L^2 ({\Bbb R})$. 
  Let the domain ${\cal D}_{\ast}\subset L^2({\Bbb R})$ consist of $y  $ such that 
     $  {\cal A}_{0}^* y \in L^2({\Bbb R})$.
The following assertion is rather standard.
  
    \begin{lemma}\label{ad}
  The operator $A_{0}$ is symmetric  on $C_{0}^\infty ({\Bbb R})$ and its adjoint $A_{0}^*$  is defined on the domain ${\cal D}(A_{0}^*) = {\cal D}_{\ast}$.     For $y \in {\cal D}(A_{0}^*)$, we have $A_{0}^* y= {\cal A}_{0}^* y $.
 \end{lemma}
 
  \begin{pf}
  By definition, ${\cal D}(A_{0}^*)$ consists of $g\in L^2 ({\Bbb R})$ such that
   \begin{equation}
  ( A_{0} g,y ) = (   g,y_{*} ) 
 \label{eq:ad1}\end{equation}
 for all $g \in C_{0}^K  ({\Bbb R})$ and some $ y_{*}\in L^2 ({\Bbb R})$; in this case $ y_{*}= A_{0}^* y$. Observe that the left-hand sides of \e{eq:ad} and  \e{eq:ad1} coincide. If ${\cal A}_{0}^* y \in L^2 ({\Bbb R})$, then \e{eq:ad1} is satisfied with $y_{*}= {\cal A}_{0}^* y$. Conversely, if  \e{eq:ad1} is satisfied, then 
 \[
  \la   g , {\cal A}_{0}^* y \ra= (   g, y_{*} ),\q  \forall  g \in C_{0}^K ({\Bbb R}),
  \]
  and hence $y_{*}= {\cal A}_{0}^* y$ so that ${\cal A}_{0}^* y \in L^2 ({\Bbb R})$.
      \end{pf}
      
      Under additional assumptions on $v(\xi)$ the operator $A_{0}^*$ is symmetric. The proof of this result requires the following auxiliary assertion. Recall that $D=id/d\xi$.
      
          \begin{lemma}\label{ad1}
  Suppose that $v \in C^1 ({\Bbb R})$, $v \in L^\infty ({\Bbb R})$, $v(\xi)>0$ and
    \begin{equation}
 | v'(\xi)| \leq C v(\xi).
 \label{eq:ad1x}\end{equation}
  Let  $z\in {\Bbb C}$ and let  $|\Im z|$ be sufficiently large. Then  the operator $v (D-z)^{-1} v^{-1}$ in $L^2 ({\Bbb R})$ defined on functions with compact supports extends to a bounded operator.
 \end{lemma}

   \begin{pf}
   Let $z=a+ib$. Since 
   \[
((D-z)^{-1}   y_{a}) (\xi)= e^{-ia\xi } ((D-i b)^{-1} y ) (\xi),  \q {\rm where} \q y_{a}(\xi)=e^{-ia\xi } y(\xi),
\]
we can suppose that $z=ib$ where $b\in {\Bbb R}$. We have to check the inequality
  \begin{equation}
 \| v (D- ib)^{-1} v^{-1} w \| \leq C \| w \|
 \label{eq:ad2}\end{equation}
 on a dense in $L^2 ({\Bbb R})$ set of elements $w $ with compact supports. Consider $w = v (D- ib)  u$ where $u \in C_{0}^\infty ({\Bbb R})$ is arbitrary; then  $w \in C_{0}^1 ({\Bbb R})$. The    set of such elements  $ w $   is dense in $L^2 ({\Bbb R})$. Indeed,
 suppose that
 \[
 (v (D- ib)  u, g_0)=0
 \]
 for some $g_{0}\in L^2 ({\Bbb R})$ and all $u\in C_{0}^\infty ({\Bbb R})$. Then $(D+ ib)  (vg_{0})=0$ and hence 
 \[
 v(\xi) g_{0}(\xi)=c e^{-b\xi}.
 \]
 Since $v \in L^\infty ({\Bbb R})$, it implies that $c=0$, and the equality $v(\xi) g_{0}(\xi)=0$ implies that $ g_{0}(\xi)=0$ because $v(\xi)  \neq 0$.

 For $w=  v (D- ib)  u$,  \e{eq:ad2} is equivalent to the inequality 
   \begin{multline}
\int_{-\infty}^\infty v^2 (\xi ) |u(\xi)|^2 d\xi \leq C
\int_{-\infty}^\infty v^2 (\xi ) |u'(\xi)- b u(\xi)|^2 d\xi
 \\
=  C\Big(
\int_{-\infty}^\infty v^2 (\xi ) |u'(\xi) |^2 d\xi- 2 b\Re\int_{-\infty}^\infty v^2 (\xi ) u'(\xi)  \bar{u}(\xi) d\xi + b^2 \int_{-\infty}^\infty v^2 (\xi ) |u(\xi)|^2 d\xi \Big).
 \label{eq:ad3}\end{multline}
 Integrating in the second term in the right-hand side by parts and using condition
  \e{eq:ad1x}, we see that 
  \[
  \Re\int_{-\infty}^\infty v^2 (\xi ) u'(\xi)  \bar{u}(\xi) d\xi = -  \int_{-\infty}^\infty v' (\xi ) v(\xi) |u(\xi)|^2 d\xi
  \]
  is bounded by
  \[
  C  \int_{-\infty}^\infty v^2 (\xi ) |u(\xi)|^2 d\xi.
  \]
  This proves inequality \e{eq:ad3} if $b$ is large enough.
   \end{pf}
   
       \begin{corollary}\label{ad1x}
       Let $b\in{\Bbb R}$ be sufficiently large. Then
  for all $k=0,1,\ldots, K$,  the operator $v D^k (Q(D) -ib)^{-1} v^{-1}$ in $L^2 ({\Bbb R})$ defined on functions with compact supports extends to a bounded operator.
 \end{corollary}
 
   \begin{pf}
   The equation $Q(x)=ib$ has $K$ solutions $x_{1}(b),\ldots, x_{K}(b)$ which for large $b$ are close to the solutions of the equation $q_{K} x^K =ib$. Therefore the roots $x_{\ell}(b) $  are simple and $|\Im x_{\ell}(b)|\to \infty $ as $b\to \infty$ for all $\ell=1,\ldots K$. Let us expand the function $x^k (Q(x) -ib)^{-1}  $ in a linear combination of the functions $(x-x_{\ell})^{-1}$ and of the constant term $1$ (for $k=K$). We can apply  Lemma~\ref{ad1} to every term $v   (D -x_{\ell}(b))^{-1} v^{-1}$. The contribution of $1$ gives the identity operator.
   \end{pf}
   
   Recall that ${\cal D}(A_{0}^*)={\cal D}_{*}$ according to  Lemma~\ref{ad}.
    Below we need additional information on this set.  Let us accept the following
   
    \begin{assumption}\label{qq}
 The function  $v \in C^K ({\Bbb R})$, $v \in L^\infty ({\Bbb R})$, $v(\xi)>0$  and estimate \e{eq:ad1x} holds.
 \end{assumption}

    \begin{lemma}\label{ad2}
    Let Assumption~\ref{qq} be satisfied.
    If $g\in {\cal D}(A_{0}^*)$, then $v D^k (vg)\in L^2 ({\Bbb R})$ for all $k= 1,\ldots, K$ and, in particular, $g\in{\sf H}^K_{\rm loc} ({\Bbb R})$.  Moreover, the coercitive estimates hold:
       \[
 \| v D^k (vg) \| \leq C \big(\| v Q(D) (vg) \| + \| g\| \big), \q k=1,\ldots, K.
\]
        \end{lemma}
    
     \begin{pf}
     By  definition of ${\cal D}_{*}$, we have $  v Q(D) (vg)\in L^2 ({\Bbb R})$ and hence
     $ w := v (Q(D)-ib) (vg)\in L^2 ({\Bbb R})$ for all $b$. Observe that
     $ v D^k (vg)=( v D^k (Q(D)-ib) ^{-1} v^{-1}) w$. Thus it remains to use Corollary~\ref{ad1x}.     
        \end{pf}
        
This lemma shows that the set $ {\cal D}_{*}\subset L^2 ({\Bbb R})$  consists of functions $g\in{\sf H}^K_{\rm loc} ({\Bbb R})$ such that  $v D^k (vg)\in L^2 ({\Bbb R})$ for all $k=1,\ldots, K$.    Now it is easy to check the following assertion.
        
          \begin{lemma}\label{ad3}
          Under Assumption~\ref{qq} the set $C_{0}^K ({\Bbb R})$ is dense in   $  {\cal D}(A_{0}^*)$  in the graph-norm $\| g \|+ \| A_{0}^* g \|$.
      \end{lemma}
    
     \begin{pf}
     Let $\varphi\in     C_{0}^\infty ({\Bbb R})$ and $\varphi(\xi)=1$ for $|\xi|\leq 1$. Set
     $\varphi_{n}(\xi)= \varphi(\xi/n)$. For an arbitrary  $ g\in {\cal D}(A_{0}^*)$, we put $g_{n}=g \varphi_{n}$. Of course $\|g -g_{n}\| \to 0$ as $n\to \infty$. Set $u =v g$, $u_{n}=v g_{n}$. We have to show that $\|v Q(D) (u- u_{n})\|\to 0$ as $n\to \infty$ or that
      \begin{equation}
 \lim_{n\to \infty} \|v D^k (u-u_{n})\| =0, \q k=0,1,\ldots, K. 
\label{eq:den}\end{equation}
Recall that $v u^{(k)} \in L^2 ({\Bbb R})$ by Lemma~\ref{ad2}. Therefore
    \[
 \lim_{n\to \infty} \|v   u^{(k)} (1-\varphi_{n})\| =0 \q {\rm and} \q  \lim_{n\to \infty} \|v   u^{(k)  }\varphi^{(l)}_{n} \| =0 
\]
for all $k=0,1,\ldots, K$ and $l\geq 1$. These relations imply \e{eq:den}.
        \end{pf}
        
        Lemma~\ref{ad3} shows that the operator $A_{0}^*$ coincides with the closure $\bar{A}_{0} $ of the operator $A_{0}$.  This yields the following assertion.

  \begin{theorem}\label{1y}
Let Assumption~\ref{qq} be  satisfied.   Then the operator $A_{0}$ defined by formula \e{eq:LOGz}  on  $C_{0}^K ({\Bbb R})$ is essentially self-adjoint. Its closure $\bar{A}_{0}=:A$ is self-adjoint on  the set ${\cal D}_{\ast}= : {\cal D} (A)$ and $Ag= v Q( D ) ( v g ) $ for $g \in{\cal D}_{\ast}$.
 \end{theorem}

    For $K$   even,   it is also possible to define $A$ in terms of the quadratic form
    \begin{equation}
 (Ag,g) =  \int_{-\infty}^\infty  \big( Q(D )( v g ) \big) v  \bar{g }d\xi .
\label{eq:LO}\end{equation}
      We suppose that  $  q_{K}>0$; then the form   $(Ag,g)  + c \| g \|^2$ is positive-definite for a sufficiently large $c>0$.
        Similarly to Theorem~\ref{1y}, it can be verified that this form defined on $C_{0}^K ({\Bbb R})$ admits the closure, and it is closed on the set $\wt{\cal D}_{\ast}$ of functions $g\in L^2({\Bbb R})$ such that  $D^k ( v g ) \in L^2({\Bbb R})$ for all $k=1,\ldots, K/2$. Then the operator $A+c I$ can be defined as a self-adjoint operator corresponding to this closed form.  Note that $\wt{\cal D}_{\ast}={\cal D} (\sqrt{A+ cI})$.


  \medskip
 
  {\bf 3.5.}
  Let us return to Hankel operators. We recall that according to Theorem~\ref{1z} the Hankel operator $H$ with kernel \e{eq:LOG} is unitarily equivalent to   differential operator  \e{eq:LOGz} where $v$ is defined by formula \e{eq:LOGz1} and $Q(x)$ is   polynomial \e{eq:LO4} with the coefficients defined by formula \e{eq:LOG5a}. To be more precise, the operators $H$ and $A$ are linked by relation  \e{eq:LOGz3}  where $F$ is operator \e{eq:MAID1}. In particular, we have
     \[
     {\cal D}(H)= F^*   {\cal D}(A) \q {\rm and}\q   {\cal D}(\sqrt{H+ c I})= F^*   {\cal D}(\sqrt{A+ c I})
     \q {\rm for}\q K \q {\rm even}.
     \] 
     Therefore the following result is a direct consequence of
     Theorem~\ref{1y}. Recall that
     the set  ${\cal D}_{0}$  consists of functions $f(t)$ satisfying estimate \e{eq:assf} for $m=0$ and some  $n> K+1$.
     
       \begin{theorem}\label{1yh}
  Let kernel $h(t)$ be defined by formulas \e{eq:LOG} and \e{eq:LOG4} where $p_{k}=\bar{p}_{k}$ for $k=0,1,\ldots, K$.  The Hankel operator $H$ with   kernel $h(t)$ is essentially self-adjoint on the domain ${\cal D}_{0}$, and its closure is self-adjoint  on the domain 
    $ F^* {\cal D}_{\ast}$.
 \end{theorem}

 \medskip
 
  \section{Spectral results}
  
  Here  we study   spectral properties of the operators $A$ and $H$.    
 
  \medskip
 
  {\bf 4.1.}
  We recall that the precise definition of the operator $A$ was given in Theorem~\ref{1y}.
  The following result relies on a construction of trial functions.

   \begin{theorem}\label{esssp}
     Let   Assumption~\ref{qq} be satisfied.
   Suppose additionally that   
     \begin{equation}
v^{(k)}  v ^{-1} \in L^\infty ({\Bbb R}), \q k=1,\ldots, K-1,
\label{eq:as}\end{equation}
 and that, for some $\d>0$,
     \begin{equation}
\Big(    \int_{-n(1+\d)}^{n(1+\d)} v(\xi)^{-2+ 4/K} d\xi \Big) \Big(   \int_{-n }^{n } v(\xi)^{-2 } d\xi   \Big)^{-1}\to 0
\label{eq:as1}\end{equation}
as $n\to \infty$.
 If $K$ is odd, then $\spec (A)={\Bbb R}$.
  If $K$ is even and $ q_{K}>0$, then $ [0,\infty) \subset \spec (A)$.
 \end{theorem}
 
   \begin{pf}
   We shall construct Weyl sequences for all $\lambda\in {\Bbb R}$ in the case of odd $K$ 
   and for all $\lambda\in [0,\infty) $ in the case of even $K$. Let $\varphi= \bar{\varphi}\in C_{0}^\infty ({\Bbb R})$, $\varphi(\xi)=1$ for $|\xi|\leq 1$ and $\varphi(\xi)=0$ for $|\xi|\geq  1+\d$. We set
    \begin{equation}
  G (\xi ; \lambda)=   \lambda^{1/K} \int_{0}^\xi v(\eta)^{-2/K} d\eta
\label{eq:LO3}\end{equation}
   and 
   \[
   g_{n} (\xi)= v (\xi)^{-1} e^{i G (\xi; \lambda)} \varphi_{n} (\xi) \q {\rm where}\q \varphi_{n} (\xi)=\varphi (\xi/n).
   \]
   Obviously,  we have
    \begin{equation}
   \| g_{n} \|^2\geq \int_{-n}^n v (\xi)^{-2} d\xi .
\label{eq:LO1}\end{equation}

Let us calculate
 \begin{equation}
   A g_{n} - (-1)^K q_{K}  \lambda g_{n}= v Q(D ) ( e^{i G  } \varphi_{n})- (-1)^K q_{K} \lambda v^{-1} e^{i G  } \varphi_{n}.
\label{eq:LO2}\end{equation}
Differentiating exponentials and using definition
  \e{eq:LO3} and conditions \e{eq:as}, we see that
 \[
  D^k  ( e^{i G (\xi; \lambda) } )=O (v(\xi)^{-2k/K}), \q k=0,1,\ldots, K-1,
\]
and
 \[
 D^K   ( e^{i G (\xi; \lambda) } )= (-1)^K \lambda v(\xi)^{-2}  e^{i G (\xi; \lambda) } + O (v(\xi)^{-2(K-1)/K}) 
\]
as $|\xi|\to \infty$. Estimating the functions $\varphi_{n}(\xi)$ and  their derivatives by constants, we find that
 \begin{equation}
  v(\xi) Q(D) ( e^{i G(\xi; \lambda) } \varphi_{n}(\xi)) = (-1)^K  q_{K}\lambda v(\xi)^{-1}  e^{i G (\xi; \lambda) } \varphi_{n} (\xi) + O (v(\xi)^{1 -2(K-1)/K}) . 
\label{eq:LO5}\end{equation}
 Substituting this expression into \e{eq:LO2}, we see that  the first term in the right-hand side of \e{eq:LO5} is cancelled with the second term in the right-hand side of \e{eq:LO2}. This yields the estimate
 \begin{equation}
 \| A g_{n} - (-1)^K q_{K}  \lambda g_{n} \|^2\leq C \int_{-n(1+\d)}^{n(1+\d)} v(\xi)^{2- 4(K-1)/K} d\xi.
\label{eq:LO6}\end{equation}
By virtue of condition \e{eq:as1}, it follows from \e{eq:LO1} and \e{eq:LO6} that
\[
\| A g_{n} - (-1)^K q_{K}  \lambda g_{n}\| \| g_{n}\|^{-1}\to 0
\]
as $n\to \infty$ so that $(-1)^K q_{K}  \lambda\in \spec (A)$.
    \end{pf}

 Let us discuss condition \e{eq:as1}. If $K=2$, it means that
 \[
 n^{-1}     \int_{-n }^{n } v(\xi)^{-2 } d\xi   \to \infty
 \]
 as $n\to\infty$. Since the integral here can be estimated from below by $n\min_{|\xi|\geq n/2} v(\xi)^{-2 } $, this condition is automatically satisfied provided $v(\xi)\to 0$ as $|\xi|\to\infty$.
 
 Let $K>2$. If 
 \[
 c| \xi|^{-\rho}\leq v(\xi) \leq C | \xi|^{-\rho}, \q 0 <c <C < \infty, \q \rho>0,
 \]
 then expression \e{eq:as1} is estimated by $C(\d) n^{-4\rho/K}$. Hence condition \e{eq:as1} is satisfied  in this case (for all $\d$). If 
 \[
 c e^{-\rho | \xi|}\leq v(\xi) \leq C e^{-\rho | \xi|}, \q 0 <c <C < \infty, \q \rho>0,
 \]
  then expression \e{eq:as1} is estimated by 
  $$C \exp\Big(2\rho n \big( (1- 2 K^{-1}) (1+\d) -1\big)\Big).$$
  This expression tends to zero if $\d<  2 (K-2)^{-1}$ so that condition \e{eq:as1} is again satisfied   for such  $\d$. On the other hand, condition \e{eq:as1} can be violated for $K>2$ if $v(\xi)$ tends to zero very rapidly (as $e^{-e^{|\xi|}}$, for example).
   
 For the next result, assumptions \e{eq:as} and \e{eq:as1} are not necessary.
 
 \begin{proposition}\label{Fxab}
  Let Assumption~\ref{qq} be satisfied.
Then $0$ is not an eigenvalue of the operator $A$.  
 \end{proposition}
 
  \begin{pf}
  Let $Ag=0$ for some $g\in {\cal D} (A)$. Put $u=vg$. Then  $u \in L^2 ({\Bbb R})$ because $v \in L^\infty ({\Bbb R})$.
  Since $v(\xi)>0$, we have $Q(D)u =0$. Denote by $x_{1}, \ldots, x_{K_{0}}\in{\Bbb C}$ different roots of the equation $Q(x)=0$. Then
     \begin{equation}
u(\xi)=\sum_{k=1}^{K_{0}} P_{k} (\xi) e^{-ix_{k}\xi}
 \label{eq:ww}\end{equation}
 for some polynomials $P_{k} (\xi)$. Observe that all exponentials do not decay at least at one of the infinities. Therefore function \e{eq:ww} does not belong to $  L^2 ({\Bbb R})$ unless all polynomials $P_{k} (\xi)$, $k=1, \ldots, K_{0}$,  are zeros.
   It follows that $u=0$ whence $g=0$.
    \end{pf}
    
     \medskip
     

  {\bf 4.2.}
  Let $K$ be even and $  q_{K} > 0$; then $A$ is semi-bounded from below and according to  \e{eq:LO}  we have 
    \begin{equation}
(A g,g)=\int_{-\infty}^\infty Q(x) | (\Phi^* ( v g)) (x)|^2 dx, \q \forall g \in {\cal D} (A).
 \label{eq:UNx}\end{equation}
Clearly, $A \geq 0$ if $Q(x) \geq 0$. On the other hand, if   $Q (x_{0}) < 0$ for some $x_{0}\in {\Bbb R}$, then $Q(x)< 0$ for some interval $\Delta$ centered at the point $x_{0}$. For every $N$,  we choose functions $\psi_{n}\in C_{0}^\infty ({\Bbb R})$ with $\supp \psi_{n}\subset \Delta$ for all $n=1,\ldots, N$  such that $\supp \psi_{n}\cap \supp \psi_m = \varnothing$ if $n\neq m$. The functions $g_{n}=v^{-1} \Phi \psi_{n}\in {\cal D} (A)$ and according to  \e{eq:UNx} the form 
 \[
 (A g,g)= \sum_{n=1}^N |\alpha_{n}|^2 \int_{-\infty}^\infty Q(x) |\psi_{n} (x)|^2 dx <0
 \]
  on all non-trivial linear combinations $g=\sum_{n=1}^N \alpha_{n} g_{n} $  of the functions   $g_{1}, \ldots, g_{N}$. This leads to the following result.

 \begin{theorem}\label{FDab}
  Let Assumption~\ref{qq} be satisfied.
Suppose that  $K$ is even and $  q_{K} >0$.
Then the     operator  $A$    is positive if   and only if
$Q (x) \geq 0$
for all $x\in {\Bbb R} $.  Moreover, if   $Q (x_{0}) < 0$ for some $x_{0}\in {\Bbb R}$, then the negative spectrum of $A$ is infinite.
 \end{theorem}

   Let  $K\geq 2 $ be even and $q_{K} >0$. According to Theorem~\ref{esssp}, $ [0,\infty) \subset \spec_{\rm ess}(A)$. Let us show that actually   we have the equality here.
 It follows from Theorem~\ref{FDab} that  for sufficiently large $\nu$   the operator $A_{\nu}= v  (Q(D)+ \nu )v \geq 0$.
 Thus we have to check that adding the   operator $\nu v^2$ does not change the essential spectrum of $A$.


 \begin{lemma}\label{unb1}
   In addition to the assumptions of Theorem~\ref{FDab} suppose 
   that $v(\xi)\to 0$ as $|\xi|\to\infty$. Then the  operator $v^2 (A + i)^{-1}$ is compact.  
     \end{lemma}
 
 \begin{pf} 
 Let a set of functions $g_{n}$ be bounded in the graph-norm $\|A g \|+ \| g \|$.    We have to check that   it is compact in the norm $ \| v^2 g \| $. Put $u_{n}= v g_{n}$.  Lemma~\ref{ad2} implies that  
    \begin{equation}
 \| v u_{n}' \| +
   \| v^{-1} u_{n} \| \leq C <\infty.
 \label{eq:UN1}\end{equation}
 We  have to show that the set $u_{n}$   is compact in the norm $ \| v  u_{n} \|$ or in $L^2 ({\Bbb R})$ because the function $v(\xi)$ is bounded. Since $v(\xi)\to 0$ as $|\xi|\to \infty$, the boundedness of the second term in   \e{eq:UN1} shows that the norms of  $u_{n}$ in $L^2(-\infty, - R)$ and $L^2( R,\infty)$ can be made arbitrary small uniformly in $n$ if $R$ is sufficiently large. Observe that $v (\xi)\geq c > 0$ on every compact interval, and hence the boundedness of the first term in   \e{eq:UN1} shows that the set $u_{n}$ is bounded in the Sobolev space ${\sf H}^1 (-R,R)$. It follows that
 this set is compact in $L^2(- R, R)$ for all $R<\infty$. 
  \end{pf} 
   
   \begin{corollary}\label{unb2}
 For an arbitrary $\nu$, we have 
 \[
 \spec_{\rm ess} (A+\nu v^2)= \spec_{\rm ess} (A).
 \]
 \end{corollary}
 
 Putting together this result with Theorem~\ref{esssp}, we obtain the following assertion.
   
      \begin{theorem}\label{esssp1} 
       In addition to the assumptions of Theorem~\ref{esssp} suppose 
   that $v(\xi)\to 0$ as $|\xi|\to\infty$.
  If $K$ is even and $ q_{K}>0$, then  
     \begin{equation}
\spec_{\rm ess}(A) = [0,\infty). 
 \label{eq:UN}\end{equation}
 \end{theorem} 
 
 We emphasize that equality \e{eq:UN} is due to the condition $v(\xi)\to 0$ as $|\xi|\to \infty$.
 If $v(\xi)=1$, then of course $\spec (A)   = \spec_{\rm ess}(A) = [P_{\rm min},\infty)$ where $P_{\rm min}=\min P(x)$ for $x\in {\Bbb R}$.
   
 \medskip

      {\bf 4.3.}
    Theorem~\ref{1z} allows us to reformulate the results of the previous subsections in terms of Hankel operators.  We recall that the precise definition of the operator $H$ was given in Theorem~\ref{1yh}. Since function \e{eq:LOGz1} satisfies Assumption~\ref{qq}  and conditions \e{eq:as}, \e{eq:as1}, the following result is a consequence of Theorems~\ref{esssp} and \ref{esssp1}.

      \begin{theorem}\label{SpB}
 Let kernel $h(t)$ be defined by formulas \e{eq:LOG} and \e{eq:LOG4} where $p_{k}=\bar{p}_{k}$ for $k=0,1,\ldots, K$.  Then:
 
 $1^0$ The point $0$ is not an eigenvalue of $H$.
 
$2^0$ If $K$ is odd, then $\spec (H)={\Bbb R}$.

$3^0$  If $K$ is even and $ p_{K}>0$, then $    \spec_{\rm ess} (H) =[0,\infty) $.
 \end{theorem}

  We emphasize that, for $K=1$, Theorem~\ref{SpB11} yields a much stronger result.
 
   Apparently the theory of Weyl-Titchmarsh-Kodaira does not apply to operators  \e{eq:LOGz} because $v (\xi)\to 0$ as $|\xi|\to \infty$. Nevertheless we conjecture that the spectrum of $A$ is absolutely continuous up to perhaps a discrete set of  eigenvalues. Moreover, we expect that the spectrum of $A$ is simple for odd $K$ and it has multiplicity $2$ for even $K$.

  Let $K$ be even and $  p_{K} > 0$; then $H$ is semi-bounded from below. Let us find conditions of the positivity of $H$.  Since the operators $H$ and $A$ are unitarily equivalent, the
  following result is a direct consequence of  Theorem~\ref{FDab}.

 \begin{theorem}\label{FDa}
 Let kernel $h(t)$ be defined by formulas \e{eq:LOG} and \e{eq:LOG4}  where $p_{k}=\bar{p}_{k}$ for $k=0,1,\ldots, K$. Suppose that $K$ is even and $  p_{K} >0$. Let   $Q(x)$ be   polynomial \e{eq:LO4} with coefficients \e{eq:LOG5a}.
Then the   Hankel operator  $H$    is positive if   and only if
$Q (x) \geq 0$
for all $x\in {\Bbb R} $.  Moreover, if   $Q (x_{0}) < 0$ for some $x_{0}\in {\Bbb R}$, then the negative spectrum of $H$ is infinite.
 \end{theorem}

   Theorem~\ref{FDa} shows that the positivity of the Hankel operator with kernel \e{eq:LOG} defined by a polynomial $P (x)$ is determined by another polynomial $Q (x)$ defined by formula \e{eq:LO4}. 
  Of course  the condition $Q (x)\geq  0$ is stronger than $P (x)\geq  0$. This follows, for example, from representation \e{eq:LOG5f}.

 In the case $K=2$,     the condition $Q (x) \geq 0$  reads as
    $ q_{1} ^2  \leq 4 q_{0} q_2$. By virtue of \e{eq:LOG6} it can be rewritten as
  \begin{equation}
  p_{1} ^2  + 2 \pi^2 p_2^2/3 \leq 4 p_{0} p_2.
\label{eq:LOG8}\end{equation}
 Obviously,  this condition  is stronger than the condition $ p_{1} ^2  \leq 4 p_{0} p_2$ guaranteeing that $h(t)\geq 0$.

The following assertion is a particular case of Theorem~\ref{FDa}.

\begin{proposition}\label{FDb}
  The    Hankel operator  $H$ with kernel 
\[
h(t)=(p_{0} + p_{1} \ln t+ p_2 \ln^2 t) t^{-1}, \q p_{0}=\bar{p}_{0},  \q p_1=\bar{p}_1,  \q p_2 >0,
\]
  is positive if  and only if     condition \e{eq:LOG8} is satisfied.  Moreover,   if $ p_{1} ^2  + 2 \pi^2 p_2^2/3 > 4 p_{0} p_2$, then the negative spectrum of $H$ is infinite.
 \end{proposition}

  \section{Hankel operators with discontinuous kernels}

  {\bf 5.1.}
   We here consider Hankel operators with singular kernels defined by formula 
\e{eq:DD}.
Hankel operators with such kernels are formally symmetric, and   we shall see later that they are essentially self-adjoint on $C_{0}^\infty ({\Bbb R}_{+})$. According to   \e{eq:H1} we have
   \begin{equation}
(Hf)(t)=\sum_{k=0}^K   (-1)^k h_{k} f^{(k)} (t_{0}-t) , \q t\in (0,t_{0}), 
\label{eq:DD1}\end{equation}
 and $(Hf)(t)=0$ for $t >t_{0}$.  Formula \e{eq:DD1} gives us the precise definition of the Hankel operator with distributional kernel \e{eq:DD}. 

Since  $L^2 (t_{0}, \infty )\subset \Ker H$, it suffices to study the restriction of $H$ on the subspace $L^2 (0,t_{0})$. It is again given by differential expression \e{eq:DD1} on $C^\infty$ functions vanishing in a neighborhood of the point $t=0$. Let us denote by $H_{0}$ the operator  \e{eq:DD1} in $L^2 (0,t_{0})$ with such domain ${\cal D} (H_{0})$. Recall that ${\sf H}^{K}  (0,t_{0})$ is the Sobolev class. Let  the set ${\sf D}_*\subset{\sf H}^{K}  (0,t_{0})$ consist of functions   satisfying the  boundary conditions
 \begin{equation}
  f(0)  = f' (0)=\cdots  = f^{(K-1)} (0 ) =0.
 \label{eq:EXEbc}\end{equation} 
The following assertion defines $H$ as a self-adjoint operator.

 \begin{lemma}\label{Dd}
The   operator $H_{0}$ is symmetric and essentially self-adjoint. Its closure $\bar{H}_{0}=:H$ is self-adjoint in $L^2 (0,t_{0})$  on the domain   ${\cal D} (H )={\sf D}_*$, and it acts by formula \e{eq:DD1}.
   \end{lemma}

  \begin{pf} 
  Let us   denote by $H_*$ differential operator \e{eq:DD1} considered as  a mapping
  $H_* : L^2 (0,t_{0}) \to C_{0}^\infty (0,t_{0})'$.  Notice that   $H_* : {\sf H}^{K}  (0,t_{0}) \to L^2 (0,t_{0}) $. 
  Integrating by parts and using that
   $h_{k}=\bar{h}_{k}$, we see that
     \begin{equation}
(H_* f, z)=-\sum_{k=1}^K   (-1)^k  h_{k} \sum_{l=0}^{k-1} f^{(k-1-l)} (t_{0}-t) \ov{ z^{(l)}(t)}\Big|_{t=0}^{t=t_{0}} + (f, H_*  z)
\label{eq:DD2}\end{equation}
 for all $f, z\in  {\sf H}^{K}  (0,t_{0})$. Observe that the non-integral terms here vanish if both functions $f$ and $z$ satisfy boundary conditions \e{eq:EXEbc}. Since $H_{0}f = H_\ast f$ for $f \in {\cal D} (H_{0})$, 
 it follows that $(H_{0}f, z)=  (f, H_{0} z)$ if  $f, z\in {\cal D} (H_{0})$.

 Let us construct the adjoint operator $H_{0}^*$. Suppose  that $(H_{0}f,z)=(f, z _{*})$ for all $f\in {\cal D} (H_{0})$ and some $z,z_{*}\in L^2 (0,t_{0})$. Choosing first $f\in C_{0}^\infty (0,t_{0})$ and using again \e{eq:DD2}, we see that  $(H_{0}f,z)=(f, H _{*} z)$ and hence
  $z_{*}= H_* z$. Since $z_{*}\in L^2 (0,t_{0})$, we find that $z\in {\sf H}^{K}  (0,t_{0})$. 
  
  For an arbitrary $f\in {\cal D} (H_{0})$, only the nonintegral terms in \e{eq:DD2} corresponding to $t=  t_0$ are equal to zero. Therefore it follows from \e{eq:DD2} that the sum of terms corresponding to $t=0$ is also zero,   that is, 
     \begin{multline*}
\sum_{k=1}^K   (-1)^k h_{k} \sum_{l=0}^{k-1} f^{(k-1-l)} (t_{0}) \ov{ z^{(l)}(0)}
\\
= -\sum_{p=0}^{K-1}   (-1)^p f^{(p)} (t_{0})  \sum_{l=0}^{K-p-1}   (-1)^l h_{p+l+1} \ov{ z^{(l)}(0)}= 0.
\end{multline*}
Since the numbers $f^{(p)} (t_{0}) $ are arbitrary, we obtain a system of $K$ equations 
 \begin{equation}
  \sum_{l=0}^{K-p-1}   (-1)^l h_{p+l+1} \ov{ z^{(l)}(0)} =0, \q p=0, \ldots, K-1, 
\label{eq:DD3}\end{equation}
for $K$ numbers $\ov{ z (0)}, \ov{ z'(0)}, \ldots, \ov{ z^{(K-1)}(0)}$.   The matrix corresponding to this system consists of elements $a_{p,l}$ parametrized by indices $p,l= 0, \ldots, K-1$. We have $a_{p,l}= (-1)^l h_{p+l+1} $ for $p+l\leq K-1$ and $a_{p,l}= 0$ for $p+l >K-1$. The determinant of this matrix is the product of skew diagonal elements $a_{p,l}$ where $p+l=K-1$ times the factor  $(-1)^{(K-1)K/2}$. Thus it  equals $  h_{K}^K$ which is not zero. Therefore it follows from \e{eq:DD3} that necessarily $z(0)=z' (0)= \cdots = z^{(K-1)} (0)=0$. It means that $ {\cal D} (H_{0}^*)\subset {\sf D}_*$ and $H_{0}^*z = H_*z$ for $z\in {\cal D} (H_{0}^*)$.

Conversely, using again \e{eq:DD2}, we see that $(H_* f, z)= (f, H_* z)$ for all $f,z\in  {\sf D}_*$. It follows that $ {\cal D} (H_{0}^*)={\sf D}_*$ and that the operator $H_{0}^*$ is symmetric. Hence the operator $H_{0}^{**}=\bar{H}_{0} $ is self-adjoint. 
 \end{pf}


We note that zero is not an eigenvalue of the operator $H$. Indeed, after the change of variables $t\mapsto t_0 -t$ the equation $(Hf )(t)=0$ reduces to the differential equation of order $K$.    Therefore the  unique solution of the equation $(Hf )(t)=0$ satisfying conditions \e{eq:EXEbc} is zero.

   \medskip

 {\bf 5.2.}
 Clearly,  $H^2$ is the differential  operator of order $2K$ defined by the formula
   \begin{equation}
(H^2f)(t)=\sum_{k,l=0}^K   (-1)^k  h_{k} h_{l} f^{(k+l)} ( t) 
\label{eq:DD3x}\end{equation}
on functions in ${\sf H}^{2K}  (0,t_{0})$ satisfying   boundary conditions \e{eq:EXEbc} and
 \[
\sum_{l=0}^K   (-1)^{ l}  h_{l} f^{(k+l)} (t_0) =0, \q k =0,1, \ldots, K-1.
 \] 
Of course the spectrum of the operator $H^2$ consists of positive eigenvalues of multiplicity not exceeding $K$ (because the differential equation $H^2f =\lambda f$ together with conditions \e{eq:EXEbc} has $K$ linearly independent solutions). These eigenvalues accumulate to $+\infty$ and their
 asymptotics  is given by the Weyl formula. However, to find the asymptotics of   eigenvalues of the operator $H$, we have to distinguish its positive and negative eigenvalues. For this reason, it is convenient to introduce an auxiliary operator $\wt{H}$ with symmetric (with respect to the point $0$) spectrum having the same asymptotics of eigenvalues as $H$.
 
   We define the operator $\wt{H}$  by the same formula \e{eq:DD1}
as $H$ but consider it on functions in  $ {\sf H}^K  (0,t_{0}/2)\oplus {\sf H}^K (t_{0}/2,t_{0})$  satisfying the boundary conditions
 \begin{equation}
  f^{(k)}(0 )  =  f^{(k)}(t_{0}/2-0 ), \q  f^{(k)}(t_{0}/2 +0 ) =    f^{(k)}(t_{0} ),  \q k=0,\ldots, K-1,  
 \label{eq:EXEbc1x}\end{equation}
for $ K  $ odd or 
  \begin{equation}
  f^{(k)}(0 )  =  f^{(k)}(t_{0}/2-0 ) =0, \q  f^{(k)}(t_{0}/2+0 )=    f^{(k)}(t_{0} )=0,  \q  k=0,\ldots, K/2-1,  
 \label{eq:EXEbc1y}\end{equation}
 for $ K $ even.
 The operator $\wt{H}$ is self-adjoint 
   in the space $L^2 (0,t_{0}/2)\oplus L^2 (t_{0}/2,t_{0})$,
and it is determined   by the matrix
  \begin{equation}
\wt{H}= \begin{pmatrix}
0 & H_{1,2}
\\
H_{2,1}    & 0
\end{pmatrix}  , \q  H_{1,2} = H_{2,1}^*,
 \label{eq:oplus1}\end{equation}
where $H_{2,1}: L^2 (0,t_{0}/2)\to  L^2 (t_{0}/2,t_{0}) $. The operator $H_{2,1}$ is  again given  by formula \e{eq:DD1} on functions in $ {\sf H}^{K} (0, t_{0}/2 )$ satisfying conditions \e{eq:EXEbc1x}  for $ K  $ odd or \e{eq:EXEbc1y}  for $ K $ even at the points $0$ and $t_{0}/2 -0$. It follows from  representation  \e{eq:oplus1} that the spectrum of the operator $\wt{H}$ is symmetric with respect to the point $0$ and consists of eigenvalues $\pm \sqrt{\mu_{n}} $ where  $\mu_{n}$  are eigenvalues of the operator $H_{2,1}^* H_{2,1}=:{\bf H}$. 

Obviously, the  operator $H_{2,1}^*$ is  again given  by formula  \e{eq:DD1} on functions in $ {\sf H}^{K} (t_{0}/2 , t_{0})$ satisfying conditions \e{eq:EXEbc1x}  for $ K  $ odd or \e{eq:EXEbc1y}  for $ K $ even at the points   $t_{0}/2 +0$ and $t_{0}$. The operator    ${\bf H} $ acts   in the space $L^2 (0,t_{0}/2)$   according to equality \e{eq:DD3x}  and its domain ${\cal D}
({\bf H})$ consists of functions $f\in {\cal D}(H_{2,1})$ such that $H_{2,1} f\in {\cal D}(H_{2,1}^*)$; in particular,  ${\cal D}({\bf H})\subset  {\sf H}^{2K} (0,t_{0}/2)$. If $K$ is odd, we have the boundary conditions $ f^{(k)}(0 )  =  f^{(k)}(t_{0}/2 ) $ for $ k=0,\ldots, 2K-1$. If $K$ is even, then equalities \e{eq:EXEbc1y} should be complemented by the boundary  conditions
 \begin{equation}
\sum_{l=K/2 - k}^K  (-1)^l h_{l}    f^{(l+k)}(0 ) = \sum_{l=K/2 - k}^K  (-1)^l h_{l}    f^{(l+k)}(t_{0}/2 )
= 0 
 \label{eq:bcc}\end{equation}
 for $ k=0,\ldots, K/2-1$.
 Note that conditions \e{eq:EXEbc1y} and \e{eq:bcc} at the point $0$ as well as at the point $t_{0}/2$ are linearly independent because $h_{K}\neq 0$.
 
 Let $\mu_{n}$ be eigenvalues  of the operator ${\bf H} $  enumerated in increasing order with multiplicities taken into account. According to  the Weyl formula we have 
  \[
\mu_{n} =   h_{K}^2 (2\pi t_{0}^{-1} n)^{2K} (1+O (n^{-1})), \q n\to \infty.
\] 
This yields  the asymptotics  of eigenvalues $\pm \sqrt{\mu_{n}} $ of the operator  $\wt{H}$. 
 
 Let us now observe that the operators $H$ and $\wt{H}$ are self-adjoint extensions of the same symmetric operator $H_{00}$ with finite deficiency indices $(2 K, 2 K)$. The operator $H_{00}$ can be defined by formula \e{eq:DD1} on $C^\infty$ functions vanishing in some neighbourhoods of the points $0$, $t_{0}/2$ and $t_{0}$.
  Therefore the operators $H$ and $\wt{H}$ have the same asymptotics of spectra. Thus we have obtained the following result.
 
   \begin{theorem}\label{powerex2}
   Let $H$ be the self-adjoint  Hankel operator  with kernel \e{eq:DD}. Then
   $  \Ker H =L^2 (t_{0}, \infty )$.
  The non-zero spectrum    of  $H$   consists of  infinite number of   eigenvalues $\lambda_{n}^{(\pm)}$ of   multiplicities not exceeding $K$ such that  $0 < \lambda_1^{(+)}\leq \lambda_2^{(+)}\leq\cdots \leq\lambda_n^{(+)} \leq\cdots$ and $0  >  \lambda_1^{(-)} \geq \lambda_2^{(-)} \geq\cdots \geq\lambda_n^{(-)} \geq\cdots$.
   Eigenvalues $\lambda_{n}^{(\pm)}$ accumulate to $\pm\infty$ as $n\to\infty$ and have  the asymptotics 
      \[
\lambda^{(\pm)}_{n} =  \pm  |h_{K}| (2\pi t_{0}^{-1} n)^K (1+O (n^{-1}))
\]
as $n\to \infty$.
The corresponding eigenfunctions $f_n^{(\pm)} (t)$ satisfy the equation
 \[
 \sum_{k=0}^K   (-1)^k h_{k} \frac{d^k f_n^{(\pm)} (t)}{dt^k}= \lambda^{(\pm)}_{n} f_n^{(\pm)} (t_{0}-t), \q t\in (0,t_{0}), 
\]
and boundary conditions \e{eq:EXEbc}.
 \end{theorem}

  \begin{remark}\label{pwerex2}
   In the case $h(t)=\d' (t-t_{0})$ we have the explicit formulas
      \[
\lambda^{(+)}_{n} = 2\pi  t_{0}^{ -1}   ( n-1/4)  ,\q 
\lambda^{(-)}_{n} = - 2 \pi  t_{0}^{ -1}   ( n-3/4) .
\]
\end{remark}


 \end{document}